\documentclass{article}
\newtheorem{lemma}{Lemma}
\newtheorem{theorem}{Theorem}
\newtheorem{corollary}{Corollary}
\def\ds{\displaystyle}
\def\aea{ & = & }
\def\nl{\\[1\baselineskip]}
\def\f{\frac}
\def\qed{\hfill \h{\rule{6pt}{6pt}}}
\newcommand\h[1]{\ensuremath{#1}}
\newcommand\beas[1]{\begin{eqnarray*} #1 \end{eqnarray*}}
\newcommand\bea[1]{\begin{eqnarray} #1 \end{eqnarray}}
\newcommand\lrs[1]{\left[ #1 \right]}
\newcommand\lrr[1]{\left( #1 \right)}
\newcommand\lrc[1]{\left\{ #1 \right\}}
\newcommand\sbt[1]{\f{\sigma(#1)}{2 #1}}
\newcommand\nck[2]{\left( \begin{array}{c} #1 \\ #2\end{array}\right)}
\newcommand\up[1]{\rule{0in}{#1 in}}
\title{On the Sum and Product of Distinct Prime Factors of an Odd Perfect Number}
\author{Anirudh Prabhu\\
West Lafayette Jr/Sr High School\\
West Lafayette, IN 47906\\
aprabhu@purdue.edu}
\date{}
\begin{document}
\maketitle
{\renewcommand{\thefootnote}{}}
\begin{abstract}
We present  lower bounds on the sum and product of the distinct
prime factors of an odd perfect number, which provide a lower
bound on the size of the odd perfect number as a function of the
number of its distinct prime factors.
\end{abstract}

\vspace{0.3in}

The study of perfect numbers\footnote{\h{N} is a {\it perfect number} if
the sum of its proper divisors (divisors less than the number) is \h{N}
itself.} dates back to Book IX of {\it The Elements} by Euclid circa 300 B.C.\ \cite{dick}.
To date, 47 perfect numbers have been discovered and all of them are even \cite{mers}.
The quest to find odd perfect numbers started more than 350 years ago, as evidenced
by the 1638 communications between Descartes and Mersenne \cite{dick},  and
has remained unsuccessful to this day; at present, the non-existence of odd
perfect numbers has not been proved either although it is known that if an odd perfect number exists it
has to exceed $10^{300}$ \cite{bren}. Over the last 350 years a significant body of work
by some of the most eminent mathematicians has focused on the conditions that must be satisfied by
odd perfect numbers. Euler showed that an odd perfect number \h{N} must have the
form \h{N=P^nQ^2} where \h{P} is a prime number and \h{P\equiv n\equiv 1} (mod 4)
and a similar result was also derived by Frenicle in 1657 \cite{dick}.
Sylvester's 1887 conjecture that an odd perfect number must have at least
six distinct prime factors was proved by Gradshtein in 1925 \cite{baco}; the result
was subsequently improved by Nielsen who showed that every odd perfect number must have
at least nine distinct prime factors \cite{niel}. In 1888 Catalan showed that an
odd perfect number \h{N} that is not divisible by 3, 5 or 7 must have at least 26
distinct prime factors; subsequently Norton improved the result by showing that
\h{N} must have at least 27 factors \cite{nort}.  Norton also showed that an odd
perfect number that is not divisible by 3 or 5 must have at least 15 distinct
prime factors \cite{nort} while Nielsen \cite{niel} showed that an odd perfect
number that is not divisible by 3 must have at least 12 distinct prime factors.
Hare showed that the number of prime factors of an odd perfect number, counting
multiplicity, is at least 75 \cite{hare}.
In 1896 Stuyvaert observed that an odd perfect number must be a sum of two squares
\cite{dick}), while Touchard \cite{touc} showed that if an odd perfect number
\h{N} exists then \h{N \equiv 1}(mod 12) or \h{N\equiv 9}(mod 36). Iannucci and Jenkins
have showed that the largest three factors of an odd perfect number must exceed
\h{10^8, 10^4} and \h{10^2} respectively \cite{ian9,ian0,jenk}. In 1913 Dickson showed that
for every positive integer \h{r}, there can only be finitely many odd perfect numbers
with \h{r} distinct prime factors \cite{dick2}. It is known that if
an odd perfect number \h{N} has \h{r} distinct prime factors then $N< 2^{4^r}$\ \cite{niel0}, and
\h{p_i < 2^{2^{i-1}}(r-i+1)}, where $2\leq i\leq 6$ and \h{p_i} is the \h{i^{th}}
smallest distinct prime factor \cite{kis1}.  Perisastri \cite{peri} showed that
the smallest prime factor $p_1$ of an odd perfect number with \h{r} distinct prime factors
must satisfy \h{p_1 \leq \f23 r + 3}.  Cohen showed that an odd perfect number must
have a factor of the form \h{p^{n} > 10^{20}}, where \h{p} is a prime number \cite{cohe2}.
Many of the results mentioned above are improvements of earlier work, the references
to which can be found in the cited papers. Details of an ongoing search for the odd perfect numbers
can be found at \cite{org}.

We derive below  lower bounds on the sum and product of the distinct prime factors
of  odd perfect numbers. If
\bea{N=\prod_{i=1}^r p_i^{n_i}\label{defn}}
is an odd perfect number, where \h{p_1, \ldots, p_r} are prime numbers, then define
\bea{\sigma(N) := \sum_{d|N} d; \hspace{0.25in} \alpha(N) := \prod_{i=1}^r p_i; \hspace{0.25in} \beta(N) := \sum_{i=1}^r p_i
 \label{defpialpha}}
\h{\sigma(N)}, the sum of all divisors of \h{N}, is called the {\it divisor}
function in the literature. We begin by proving the following lemma.
\begin{lemma}  If \h{N>1} is an odd perfect number then \h{\sigma(\alpha(N)) < 2\ \alpha(N)}.
\end{lemma}
{\bf Proof:}  Let \h{M:=\alpha(N)}. Clearly $M\leq N$ and first we show that $M<N$.  If \h{M=N}, then
by Euler's theorem \h{N=M=A^n Q^2} where \h{A} is a prime
number and \h{A\equiv n \equiv 1} (mod 4). Since \h{M} is a
product of distinct primes, the uniqueness of prime factorization implies
that it cannot contain a factor that
is a perfect square, and therefore $N=M=A^n$.  Since \h{A} is
prime and every prime occurring in \h{M} has a unit exponent,
\h{n=1}, showing that \h{N} must be a prime number.  But if
\h{N} is a prime number, \h{\sigma(N) = N+1 \neq 2N} for
$N>1$, contradicting the assumption that \h{N} is a perfect
number.  Therefore we conclude that $M<N$.  The above argument
also proves that if \h{N} is an odd perfect number, as described
in the Lemma, then the
condition $n_1 = n_2 = \ldots = n_r=1$ cannot hold, a conclusion
that is also implied by a result of Steuerwald's \cite{steu}.

Next consider an arbitrary number \h{B} whose prime
factorization is of the form
\bea{B = p\cdot \prod_{i=1}^s q_i^{m_i} \label{defb}}
where \h{p, q_1, \ldots, q_s> 1} are prime numbers.  For
$n>1$, let
\bea{C = p^n\cdot \prod_{i=1}^s q_i^{m_i}. \label{defc}}
Then
\beas{
\ds \f{\sigma(C)}{2C} \aea \ds \f{\left(\sum_{j=0}^n p^j \right)\prod_{i=1}^s \left( \sum_{k=0}^{m_i} q_i^k\right)}{2 p^n \cdot
\prod_{i=1}^s q_i^{m_i}}
= \ds \left[\f{\left(\sum_{j=0}^1 p^j \right)\prod_{i=1}^s \left( \sum_{k=0}^{m_i} q_i^k\right)}{2 p \cdot
\prod_{i=1}^s q_i^{m_i}} \right]
\cdot \left[
\ds \f{\sum_{l=0}^n p^l}{\sum_{k=0}^1 p^k}
\right]\left[
\ds \f{1}{p^{n-1}}\right]\nl
\aea \left[ \ds \f{\sigma(B)}{2B}\right]
\left[
\ds \left(\f{p^{n+1}-1}{p-1}\right) \left(\f1{1+p}\right)
\right] \f1{p^{n-1}}\nl
\aea \left[ \ds \f{\sigma(B)}{2B}\right] \left[
\f{p^{n+1} - 1}{p^{n+1} - p^{n-1}}
\right]
}
Since \h{p,n >1}, we have \h{p^{n+1} - 1 > p^{n+1} - p^{n-1}}
and therefore we conclude that for numbers \h{B} and \h{C}
defined as in (\ref{defb}) and (\ref{defc}),
\bea{
\ds\f{\sigma(C)}{2C} > \f{\sigma(B)}{2B} \label{keyineq}
}
If we consider the sequence \h{M=B_0, B_1, \ldots, B_r=N}, where for \h{1\leq k\leq r}
\bea{
B_k = \prod_{i=1}^k p_i^{n_i} \prod_{j=k+1}^r p_j
} then repeated application of inequality (\ref{keyineq}) shows that
\bea{
\sbt{M} <  \sbt{B_1} < \sbt{B_2} < \ldots < \sbt{B_r} = \sbt{N} \label{mln}
}
Since $M\neq N$, and \h{N} being a perfect number $\sigma(N) = 2N$, using (\ref{mln}) we conclude that
\\

\h{\hspace{1.8in}
\ds \f{\sigma(M)}{2M} < \f{\sigma(N)}{2N} = 1\hspace{0in}
\hfill \h{\rule{6pt}{6pt}}}

\ \\

The following theorem is the main result.
\begin{theorem}
If \h{N>1} is an odd perfect number with \h{r} distinct prime factors then
\beas{\alpha(N) > \ds \f{1}{\rule{0in}{0.2in}\lrr{  2^{^{ \f 1 r}} - 1}^r}; \hspace{0.25in} \beta(N) > \f{r}{\rule{0in}{0.15in}2^{\f 1 r} - 1}}
\end{theorem}
{\bf Proof:} Let \h{N} be defined as in (\ref{defn}) with distinct prime factors
\h{p_1, \ldots, p_r}.  For \h{1\leq k\leq r}, define
\bea{
S_k  :=   \ds \sum_{1 \leq i_1 < i_2 < \ldots < i_k \leq r}  \f1{\left(
p_{i_1} \cdot p_{i_2} \cdot \ldots \cdot p_{i_k}
\right)} \label{sumsk}
}
The sum in (\ref{sumsk}) is over all \h{k}-subsets of \h{\{p_1, \ldots, p_r\}}.  Since
each of the \h{p_1, \ldots, p_r} occurs in exactly \h{\nck{r-1}{k-1}} terms
in the sum, using the GM-HM inequality we have
\beas{
\ds \f{\nck r k}{S_k} < \lrs{\lrr{\prod_{i=1}^r p_i}^{\nck{r-1}{k-1}}}^{\f1{\nck r k}} = \lrs{\rule{0in}{0.2in}\alpha(N)}^{\ds\f k r}
}
The GM-HM inequality is strict since we are considering distinct prime numbers. The
above inequality can be rewritten as
\bea{S_k > \nck r k  \lrs{\rule{0in}{0.2in}\alpha(N)}^{-\ds\f k r}\label{insk}}
Using Lemma 1 and (\ref{insk}) we get
\beas{
1 > \ds \f{\sigma(\alpha(N))}{2\ \alpha(N)} \aea \ds \f{\ds\prod_{i=1}^r (1+p_i)_{\rule{0in}{0.05in}}}{\rule{0in}{0.2in}2\cdot \ds\prod_{i=1}^r p_i}  = \f12 \lrc{1 + \sum_{k=1}^r S_k}
> \f12 \lrc{1 + \sum_{k=1}^r \nck r k \lrs{\lrr{\alpha(N)}^{-\ds\f1r}}^k}=
\f12 \lrc{1+\lrr{\alpha(N)}^{-\ds\f1r}}^r
}
which implies that
\beas{
1 > \f12 \lrc{1+\lrr{\alpha(N)}^{-\ds\f1r}}^r
}
or
\bea{
\alpha(N) > \ds\f1{\up{0.2}\lrr{2^{\f 1 r}-1}^r} \label{inab}
}
as claimed.

The bound on \h{\beta(N)} can be derived using the AM-GM inequality,
applying which  and using (\ref{inab}) we get
\beas{
\ds \f{\beta(N)}{r} = \f{\ds\sum_{i=1}^r p_i}{r} > \lrc{\ds\prod_{i=1}^r p_i}^{\f1r} = \lrc{\alpha(N)}^{\f1r}
> \ds\f1{\lrr{2^{\f 1 r}-1}}
}
from which we obtain the inequality
\bea{
\beta(N) > \f{r}{\lrr{2^{\f 1 r}-1}}
}
as claimed. \qed

For an odd perfect number \h{N}, since $N > \alpha(N)$  we immediately have the
corollary
\begin{corollary}
If \h{N} is an odd perfect number with \h{r} distinct prime factors then
\beas{N>\ds\f1{\up{0.2}\lrr{2^{\f1r}-1}^r}}
\end{corollary}
Together with Nielsen's result \cite{niel0} the current lower and upper bounds
on an odd perfect number with \h{r} distinct prime factors can be summarized
as
\beas{
\ds\f1{\up{0.2}\lrr{2^{\f1r}-1}^r} < N < \ds 2^{4^r}
}

Lemma 1  also yields a simple derivation of an upper bound on the sum of
reciprocals of the distinct prime factors of an odd perfect number.  The upper
bounds presented below are slightly weaker than those reported in \cite{cohe,suha}
but the derivations are considerably shorter.
\begin{theorem}
If \h{N>1} is an odd perfect number with prime factorization
\h{\ds N = \prod_{i=1}^r p_i^{n_i}}, where \h{p_1, \ldots, p_r}
are distinct prime numbers, then
\beas{
\f1{p_1} + \f1{p_2} + \ldots + \f1{p_r} < 1
}
\end{theorem}
{\bf Proof:} Let \h{M:=\alpha(N)}.
From Lemma 1, we have
\bea{
\ds\frac{\sigma(M)}{2M} \aea \ds\frac{\prod_{i=1}^r \left( 1 + p_i \right)}{2\prod_{i=1}^r p_i} = \f12 + \f12 \left(
\sum_{i=1}^r \f1{p_i}\right) + \f12 \left(
\sum_{1\leq i < j \leq r} \f1{p_i p_j}\right) + \ldots +
\f12 \left( \f1{p_1 \ldots p_r}\right) < 1 \label{mainineq}
}
Since \h{p_1, \ldots, p_r > 0},
\beas{\f12 + \f12 \left(
\sum_{i=1}^r \f1{p_i}\right) < 1}
from which the claim in the theorem follows. \hfill \h{\rule{6pt}{6pt}}

\vspace{12pt}

The bound in Theorem 1 can be improved as follows.
\begin{theorem}
If \h{N>1} is an odd perfect number that has the prime factorization
\h{\ds N = \prod_{i=1}^r p_i^{n_i}}, where\\ \h{p_1< p_2< \ldots< p_r=P}
are distinct prime numbers, then
\beas{
\f1{p_1} + \f1{p_2} + \ldots + \f1{p_r} < 1 - \lrc{\lrs{1+\f1 P}^r - \lrs{1+\f r P}}
}
\end{theorem}
{\bf Proof:}  Let {\bf a} = \h{\{a_1, \ldots, a_r\}} be a set of \h{r}
distinct positive real numbers. Without loss of generality we will assume that
\h{0<a_1 < a_2 < \ldots < a_r}.  For \h{1\leq k\leq r}, define
\bea{
S_k({\bf a}) =   \ds \sum_{1 \leq i_1 < i_2 < \ldots < i_k \leq r}  \left(
a_{i_1} \cdot a_{i_2} \cdot \ldots \cdot a_{i_k}
\right) \label{defsk}
}
That is, \h{S_k({\bf a})} is the sum of the products of subsets of \h{k} numbers chosen
from the \h{r} numbers in {\bf a}. Observe that each number, such as \h{a_1}, occurs
in \h{\left( \begin{array}{c} r-1\\ k-1 \end{array}\right)} products in the sum (\ref{defsk}).
Therefore, using the AM-GM inequality we have
\beas{
\f{S_k({\bf a})}{\left( \begin{array}{c} r\\ k \end{array}\right)} & \geq &
\lrr{\rule{0in}{0.2in} a_1\cdot \ldots \cdot a_r}^{\f{\left( \begin{array}{c} r-1\\ k-1 \end{array}\right)}{\left( \begin{array}{c} r\\ k \end{array}\right)}} = \lrr{\rule{0in}{0.15in} a_1\cdot \ldots \cdot a_r}^{\ds \f k r} >
\lrr{a_1}^k \hspace{0.2in} \Rightarrow \hspace{0.2in} S_k({\bf a}) \geq \left( \begin{array}{c} r\\ k \end{array}\right) \lrr{a_1}^k
}
using which we get
\bea{
\sum_{k=2}^r S_k({\bf a}) \geq \sum_{k=2}^r \left( \begin{array}{c} r\\ k \end{array}\right) \lrr{a_1}^k = \lrr{1+a_1}^r
- (1 + ra_1) \label{bnd}
}
Let \h{M:=\alpha(N)}
and set \h{\ds a_1 = \f1{p_r}=\f1P, \ldots a_r=\f1{p_1}}. Then from (\ref{mainineq}), we have
\bea{
\ds\frac{\sigma(M)}{2M} \aea \ds\frac{\prod_{i=1}^r \left( 1 + p_i \right)}{2\prod_{i=1}^r p_i} = \f12 + \f12 \left(
\sum_{i=1}^r \f1{p_i}\right) + \f12 \left(
\sum_{1\leq i < j \leq r} \f1{p_i p_j}\right) + \ldots +
\f12 \left( \f1{p_1 \ldots p_r}\right)  \nonumber \nl
\aea \f12 + \f12 S_1({\bf a}) + \f12 S_2({\bf a}) + \ldots +
\f12 S_r({\bf a})\nonumber \nl
& < & 1 \label{fbnd}
}
Using (\ref{bnd}) and (\ref{fbnd}) we obtain
\beas{
\f1{p_1} + \f1{p_2} + \ldots + \f1{p_r} < 1 - \lrc{\lrs{1+\f1 P}^r - \lrs{1+\f r P}}
}
as claimed \hfill \h{\rule{6pt}{6pt}}

\end{document}